\documentclass[11pt]{amsart}
\usepackage{amssymb,amsmath}

\newcommand{\C}{{\mathbb C}}
\newcommand{\Q}{{\mathbb Q}}

\newcommand{\Z}{{\mathbb Z}}
\newcommand{\N}{{\mathbb N}}

\newcommand{\OO}{{\mathcal O}}

\newcommand{\pp}{{\mathfrak p}}
\newcommand{\qq}{{\mathfrak q}}

\newcommand{\INT}{\operatorname{INT}}

\newcommand{\F}{{\mathbb F}}

\newcommand{\ord}{{{\operatorname{ord}}}}

\newcommand{\sep}{{{\operatorname{sep}}}}
\newcommand{\perf}{{{\operatorname{perf}}}}

\newcommand{\isom}{\cong}
\newcommand{\tensor}{\otimes}

\newcommand{\tr}{{\operatorname{tr}}}
\newcommand{\nr}{{\operatorname{nr}}}
\newcommand{\Br}{\operatorname{Br}}

\newcommand{\comment}[1]{}
\newtheorem{thm}{Theorem}[section]
\newtheorem{lem}[thm]{Lemma}

\newtheorem{prop}[thm]{Proposition}

\theoremstyle{definition}
 
\newtheorem{defn}{Definition}
\theoremstyle{remark}

\begin{document}
\title[Integrality at a Prime]{Integrality at a Prime for Global
  Fields and the Perfect Closure of Global Fields of Characteristic
  $p>2$}
\author{Kirsten Eisentr\"ager}
\address{School of Mathematics, Institute for Advanced
  Study, Einstein Drive, Princeton, NJ 08540}
\email{eisentra@ias.edu}
\begin{abstract}
  
  Let $k$ be a global field and $\pp$ any nonarchimedean prime of $k$.
  We give a new and uniform proof of the well known fact that the set
  of all elements of $k$ which are integral at $\pp$ is diophantine
  over $k$.  Let $k^{\perf}$ be the perfect closure of a global field
  of characteristic $p>2$. We also prove that the set of all elements
  of $k^{\perf}$ which are integral at some prime $\qq$ of $k^{\perf}$
  is diophantine over $k^{\perf}$, and this is the first such result
  for a field which is not finitely generated over its constant field.
  This is related to Hilbert's Tenth Problem because for global fields
  $k$ of positive characteristic, giving a diophantine definition of
  the set of elements that are integral at a prime is one of two steps
  needed to prove that Hilbert's Tenth Problem for $k$ is undecidable.
\end{abstract}

\thanks{The research for this paper
 was done while the author was at the University of California,
 Berkeley.}  
\maketitle
\section{Introduction}\label{theintro}
Hilbert's Tenth Problem in its original form was to find an algorithm
to decide, given a polynomial equation $f(x_1,\dots,x_n)=0$ with
coefficients in the ring $\Z$ of integers, whether it has a solution
with $x_1,\dots,x_n \in \Z$.  
Matijasevi{\v{c}} \cite{Mat70}, building on earlier work by Davis,
Putnam, and Robinson \cite{DPR61}, proved that no such algorithm
exists, {\it i.e.}\ Hilbert's Tenth Problem is undecidable.

Since then, analogues of this problem have been studied by asking the
same question for polynomial equations with coefficients and solutions
in other commutative rings $R$.  We refer to this as Hilbert's
Tenth Problem over $R$.
%
%
Perhaps the most important unsolved question in this area is Hilbert's
Tenth Problem over the field of rational numbers.  Diophantine
undecidability has been proved for several function fields of
characteristic~$0$: In \cite{Den78} Denef proves the undecidability of
Hilbert's Tenth Problem for rational function fields over formally
real fields. In 1992 Kim and Roush \cite{KR92} showed that the problem
is undecidable for the purely transcendental function field
$\C(t_1,t_2)$, and in \cite{Eis} this is generalized to finite
extensions of $\C(t_1,\dots,t_n)$ for $n \geq 2$.

Hilbert's Tenth Problem for the function field $k$ of a curve over a
finite field is also undecidable.  This was proved by Pheidas for
$k=\F_q(t)$ with $q$ odd, and by Videla \cite{Vi94} for $\F_q(t)$ with
$q$ even. In \cite{Sh96,Sh2000} Shlapentokh generalized Pheidas'
result to finite extensions of $\F_q(t)$ with $q$ odd and to certain
function fields over possibly infinite constant fields of odd
characteristic, and the remaining cases in characteristic $2$ are
treated in \cite{Eis2002}.  Before we can state the results of this
paper we need the following definition.

\begin{defn}
  1. If $R$ is a commutative ring, a {\em diophantine equation over
    $R$} is an equation $P(x_1,\dots,x_n)=0$ where $P$ is a polynomial
  in the variables $x_1, \dots, x_n$ with coefficients in
  $R$.\\
  2. A subset $S$ of $R^k$ is {\em diophantine over $R$} if there is a
  polynomial\linebreak$P(x_1,\dots,x_k,$ $y_1,\dots,y_m) \in
  R[x_1,\dots,x_k,y_1,\dots,y_m]$ such that \[S=\{(x_1, \dots, x_k)
  \in R^k: \exists \;y_1, \dots,y_m \in R,\;
  (P(x_1,\dots,x_k,y_1,\dots,y_m)=0)\}.\]
\end{defn}
When $R$ is not a finitely generated algebra over $\Z$, we restrict
our attention to diophantine equations whose coefficients are in a
finitely generated algebra over $\Z$.

For global fields of positive characteristic, Proposition~\ref{below}
below \cite[p.~319]{Sh96} is used to prove undecidability of Hilbert's
Tenth Problem.  For the purposes of this paper, global fields are
algebraic number fields or finite extensions of the rational function
fields $\F_q(t)$.
A prime of a global field $k$ is an equivalence class of nontrivial
absolute values of $k$. A nonarchimedean prime is an equivalence class
of nontrivial nonarchimedean absolute values of $k$. For a
nonarchimedean prime $\pp$ of a global field $k$ we denote by
$\ord_{\pp}$ the associated normalized additive discrete valuation
$\ord_{\pp}:k^{*} \twoheadrightarrow \Z$.
\begin{prop}\label{below}
  Let $k$ be a global field of positive characteristic, let $p$ be a
  rational prime, and let $\pp$ be a prime of $k$.  Assume that the
  sets $p(k):=\{(x,w) \in k^2: \exists s \in \N, w=x^{p^s}\}$ and
  $\INT(\pp):=\{x \in k:\ord_{\pp} x \geq 0\}$ are diophantine. Then
  Hilbert's Tenth Problem for $k$ is undecidable.
\end{prop}
So for global fields of positive characteristic, a diophantine
definition of the set of elements which are integral at some prime
$\pp$ is one of two main steps used to prove undecidability of
Hilbert's Tenth Problem.
 \comment{ Let $F$ be a
  field. The proofs of the undecidability results for fields of
  characteristic $p>0$ mentioned above have two main steps.  One of
  the steps is showing that the set of all elements which are integral
  at some fixed prime $\pp$ of $F$ is diophantine.  So it is of
  interest to study for which fields we can give a diophantine
  definition of all elements which are integral at some prime $\pp$ of
  the field.}  \comment{
  For the undecidability results for function fields $k$ of curves
  over finite fields \cite{Phei91,Vi94,Sh96,Eis2002}, all the proofs
  have two main steps. One of the steps is showing that the set of
  elements which are integral at some fixed prime $\pp$ of $k$ is
  diophantine. So it is of interest to study for which fields we can
  give a diophantine definition of all elements which are integral at
  some prime $\pp$ of the field.  }

In this paper we will prove two results. We give a different and more
uniform proof of the known fact that for any global field $k$ and any
nonarchimedean prime $\pp$ of $k$ the set of elements of $k$ which are
integral at $\pp$ is diophantine.
For number fields the result was already
implicit in the work of Robinson \cite{Rob49,Rob59}, and explicitly
written down in \cite[Proposition~3.1]{KR92b}. Their proof relies on
the Hasse principle for quadratic forms. For global function fields
the result was proved in \cite{Sh94}. There is also another approach
by Rumely~\cite{Rum80} that uses the Hasse norm principle.  Our
approach uses the Brauer group of $k$. We also prove the following new
result:
\begin{thm}\label{perfect}
  Let $k$ be a global field of characteristic $p>2$, and let
  $k^{\perf}$ be the perfect closure of $k$. Let $\pp$ be a prime of
  $k^{\perf}$.  The set $\{ x \in k^{\perf}: \ord_{\pp} x \geq 0 \}$
  is diophantine over $k^{\perf}$.
\end{thm}
The perfect closure of a field $k$ of characteristic $p$ is obtained
by adjoining $p^n$-th roots of all elements of $k$ for all $n \geq 1$.
A prime $\pp$ of $k^{\perf}$ is an equivalence class of nontrivial
absolute values of $k^{\perf}$. The associated additive valuation
$\ord_{\pp}$ is no longer discrete since every element of
$k^{\perf}$ is a $p$-th power.

The perfect closure of $\F_q(t)$ is
$K:=\F_q(t,t^{1/p},$ $t^{1/{p^2}},$ $t^{1/{p^3}},\dots)$.  We will
first prove Theorem~\ref{perfect} for $K$. Let $k$ be any global field
of characteristic $p>0$. Then $k$ is a finite extension of $\F_q(t)$
for some $q=p^n$. We will show in Section~\ref{sectionperf} that the
perfect closure $k^{\perf}$ of $k$ is also obtained by adjoining
$p^n$-th roots of $t$, and that the proof for $K$ generalizes to
$k^{\perf}$. These perfect closures are not finitely generated over
their constant fields.  This distinguishes them from all the function
fields mentioned above.

\section{Background}
In this section we will state some of the definitions and theorem
about division algebras and Brauer groups that are needed in the next
two sections.
\begin{defn}[Quaternion Algebras]
  Let $F$ be a field of characte\-ri\-stic $\neq~2$. For $a, b \in F^*$,
  let $H(a,b)$ be the $F$-algebra with basis $1,i,j,k$ (as an
  $F$-vector space) and with multiplication rules
\[
i^2 =a,j^2 = b, ij=k=-ji.
\]
Then $H(a,b)$ is an $F$-algebra which is called a {\em quaternion algebra}
over $F$.
\end{defn}
One can show that $H(a,b)$ is either a division algebra or isomorphic
to $M_2(F)$. (Here $M_2(F)$ is the algebra of $2 \times 2$ matrices.)
\begin{defn}
  1. An algebra $A$ is said to be {\em central simple} over a field
  $F$ if $A$ is a simple algebra having $F$ as its center.\\
  2. The matrix algebra $M_n(F)$ is called a {\em split} central
  simple algebra over $F$. If $A$ is a finite dimensional central
  simple algebra over $F$, then an extension field $E$ of $F$ is
  called a {\it splitting field} for $A$ if $A\tensor_F E \isom
  M_n(E)$ for some
  $n$.
\end{defn}
\begin{prop}\label{quaternion}
Let $F$ be a field of characteristic $\neq 2$. Every $4$-di\-men\-sio\-nal
central simple algebra over $F$ is isomorphic to $H(a,b)$ for some
$a,b \in F^{*}$.
\end{prop}
\begin{proof}
This is Proposition~1 in \cite[p.~128]{Bou}.
\end{proof}
In characteristic $2$ something similar holds:
\begin{prop}
\label{bourbaki}
  Let $F$ be a field of characteristic $2$.  Let $D$ be a central
  division algebra over $F$ such that for each $x \in D$, we have
  $[F(x):F] \leq 2$. Then $D$ admits a basis (1,u,v,w) over $F$ such
  that
\[
u^2 = a, v^2 = v + b, uv=w, vu = w+u , w^2 = ab, vw = bu
\] 
\[
wv = bu +w, wu = a + av, uw = av,
\]
where $a,b \in F$. We will denote this algebra again by $H(a,b)$.
\end{prop}
\begin{proof}
This is Exercise 4 in \cite[p.~130]{Bou}.
\end{proof}
\begin{defn}
Let $k$ be a global field. Let $\pp$ be a prime of $k$, and let $k_{\pp}$ be
the completion of $k$ at $\pp$. A quaternion algebra $A$ over $k$ is
said to {\em split} at $\pp$ if
\[
A \tensor _k k_{\pp} \isom M_2(k_{\pp}) \mbox{ as }k_{\pp} \mbox{-algebras}.
\]
Otherwise $A$ is {\it ramified} at $\pp$.
\end{defn}
{\bf Notation: }For any field $F$, let $F^{\sep}$ denote a separable
closure of $F$.\\
We have the following Proposition.
\begin{prop}
  Let $A$ be a finite dimensional central simple algebra over a field
  $F$. There exists an $F^{\sep}$-algebra isomorphism $\iota: A\tensor
  _F F^{\sep} \to M_r({F^{\sep}}).$ The characteristic polynomial
  $P_a(x) \in F^{\sep}[x]$ of $\iota(a \tensor 1)$ is independent of
  the choice of $\iota$. Moreover, $P_a(x) \in F[x]$.
\end{prop}
\begin{proof}
This is proved in \cite[pp.\ 113--114]{Reiner}.
\end{proof}
\begin{defn}
  Let $A$ be as above.  The {\em reduced trace} $tr(\alpha)$ of
  $\alpha \in A$ is defined to be the trace of $\iota(\alpha \tensor
  1)$, for any choice of $\iota$ as above. Similarly the {\em reduced
    norm} $nr(\alpha)$ is defined to be the determinant.
\end{defn}
We can compute the following:
\begin{lem}
  Let $H(a,b)$ be a quaternion algebra over a field $F$ of
  characteristic $\neq 2$. The reduced trace
  $\tr(x_1+x_2i+x_3j+x_4k)$ equals $2x_1$, and the reduced norm
  $\nr(x_1+x_2i+x_3j+x_4k)$ equals $x_1^2 - ax_2^2 -bx_3^2 +abx_4^2$
  for any $x_1,\dots, x_4 \in F$.
\end{lem}
\begin{lem}
  Let $D$ be a $4$-dimensional division algebra over a field $F$ of
  characteristic $2$, so that $D=H(a,b)$ as in
  Proposition~\ref{bourbaki} for some $a, b \in F^*$. Let $(1,u,v,uv)$
  be a basis of $D$ over $F$ as in Proposition~\ref{bourbaki}. For an
  element $x_1 + x_2u+x_3v+x_4uv$ we have
  $\tr(x_1+x_2u+x_3v+x_4uv)=x_3$ and $\nr(x_1+x_2u+x_3v+x_4uv)= x_1^2
  + x_1x_3 + b x_3^2 + a(x_2^2 + x_2x_4 + bx_4^2)$.
\end{lem}
\begin{proof}
  This follows from Proposition~10 in \cite[p.~144]{Bou} and from
  Exercise~6 in \cite[p.~147]{Bou}.
\end{proof}
\begin{defn}[Brauer group]
  Let $A$ and $B$ be finite dimensional central simple algebras over a
  field $F$. We say that $A$ and $B$ are similar, $A \sim B$, if $A
  \tensor_F M_n(F) \isom B\tensor_F M_m(F)$ for some $m$ and $n$.
  Define the {\em Brauer group} of $F$, $\Br(F)$, to be the set of
  similarity classes of central simple algebras over $F$, and write
  $[A]$ for the similarity class of $A$. For classes $[A]$ and $[B]$,
  define
\[
[A][B]:= [A\tensor _F B].
\]
This is well defined and makes $\Br(F)$ into an abelian group.
\end{defn}
Each similarity class of $\Br(F)$ is represented by a central division
algebra, and two central division algebras representing the same
similarity class are isomorphic \cite[p.~100]{Mi97}.
\begin{thm}\label{localfield}
Let $K$ be a nonarchimedean local field.
\begin{enumerate}
\item The Brauer group of $K$ is isomorphic to $\Q /\Z.$
\item Let $D/K$ be a division algebra of degree $n^2$. The order of
  $[D]$ in $\Br(K)$ is $n$.
\end{enumerate}
\end{thm}
\begin{proof}$\left.\right.$
\begin{enumerate}
\item This is Theorem~9.22 in \cite{Jac89}.
\item This is Theorem~9.23 in \cite{Jac89}.
\end{enumerate}
\end{proof}
\begin{thm}
\label{exact}
  Let $k$ be a global field.  There is an exact sequence
\[
0 \to \Br(k) \to \bigoplus_{ v \in M_{k}} \Br(k_v)
\to \Q /\Z \to 0,
\]
where
$M_{k}$ denotes the set of nonequivalent nontrivial absolute values of
$k$.
\end{thm}
\begin{proof}
This is Remark (ii) in \cite[p.~277]{Reiner}.
\end{proof}
\begin{prop}
  Let $K$ be a nonarchimedean local field, and let $D$ be a
  finite dimensional central division algebra over $K$. The valuation
  on $K$ has a unique extension to $D$.
\end{prop}
\begin{proof}
This is proved in \cite[p.~182]{SeLocal}.
\end{proof}
\section{Integrality at a prime for global fields}
In this section we will prove the following
\begin{thm}\label{orderproof}
  Let $k$ be a global field. Let $\pp$ be a nonarchimedean prime of
  $k$. The set $\{ x \in k :\ord_{\pp} x \geq 0\}$ is diophantine
  over $k$.
\end{thm}
\begin{proof}
We will first prove this when the characteristic of $k$ is not $2$ and
then say how the proof has to be modified in characteristic $2$.

For any nonarchimedean prime $\pp$ of $k$ let $R_{\pp}:= \{x \in k: \ord_{\pp} x \geq 0\}$.\\
\noindent{\bf Claim:} Given two distinct nonarchimedean primes $\pp$
and $\qq$ of $k$ there exists a subset $S \subseteq R_{\pp} \cap
R_{\qq}$ containing a
subgroup $G$ of finite index in $R_{\pp} \cap R_{\qq}$, such that $S$
is diophantine over $k$.\\
{\bf Proof of Claim:} By the approximation theorem we may choose $p,q
\in k$ such that $\ord_{\pp}p =1$, $\ord_{\qq}p =0$, $\ord_{\pp} q
=0$, and $\ord_{\qq} q = 1$. By Theorem~\ref{localfield} and
Theorem~\ref{exact} we can find a central division algebra $H$ that is
ramified exactly at $\pp$ and $\qq$ and which has degree $4$ over $k$.
By Proposition~\ref{quaternion}, $H\isom H(a,b)$ for some $a,b \in
k^{*}$. Let $\OO_{\pp}$ be the valuation ring of $k_{\pp}$, where
$k_{\pp}$ is the completion of $k$ at the prime $\pp$.  Let $A_{\pp}$
be the valuation ring of $H_{\pp}:= H \tensor k_{\pp}$. Then $A_{\pp}$
is a free $\OO_{\pp}$-module of rank $4$. Since $H(a,b) \isom
H(ax^2,by^2)$ for $x,y \in k^*$, we can choose $i,j \in H$ that are
integral at $\pp$ and $\qq$, and then
\[
p^rA_{\pp} \subseteq \OO_{\pp} + \OO_{\pp}i + \OO_{\pp}j +
\OO_{\pp}ij, \mbox{ and }\] 
\[
q^rA_{\qq} \subseteq \OO_{\qq} + \OO_{\qq}i + \OO_{\qq}j +
\OO_{\qq}ij \mbox{ for some }r \geq 0.
\]
Now let
\[
T:= \{x_1 \in k: (\exists \, x_2,x_3,x_4 \in k):(x_1^2 -ax_2^2 -bx_3^2
+abx_4^2 = pq)\}.
\]
Then $S=(pq)^{r}T$ has the desired property.  Suppose $x_1 \in T$.
Then there exists $\alpha = x_1 + x_2 i + x_3j + x_4 ij \in H$ whose
reduced norm equals $pq$. Since $pq \in \OO_{\pp}$ it follows that
$\alpha \in A_{\pp}$. Then $p^rx_1 \in \OO_{\pp}$.  Similarly we can
show that $q^rx_1\in \OO_{\qq}$, so $(pq)^r x_1 \in \OO_{\pp} \cap
\OO_{\qq}$. Hence $S \subseteq \OO_{\pp}\cap \OO_{\qq} \cap k = R_{\pp} \cap
R_{\qq}$.

Conversely assume that $x_1 \in k$ and that $x_1 \in pR_{\pp}
\cap qR_{\qq}$. Then the equation
\[
X^2 - 2x_1 X + pq = 0
\]
is Eisenstein at $\pp$ and $\qq$, so a root $\beta$ generates a quadratic
field extension, and $\beta$ also generates a quadratic extension
$k_{\pp}(\beta)$ of $k_{\pp}$ and a quadratic extension $k_{\qq}(\beta)$ of
$k_{\qq}$. By \cite[Remark~4.4, p.~110]{Mi97} any quadratic extension
field of the local field $k_{\pp}$ is a splitting field for $H$ over
$k_{\pp}$. Hence $k_{\pp}(\beta)$ splits $H$ locally, and by
Theorem~\ref{exact} it follows that $k(\beta)$ splits $H$.  Since $k(\beta)$
splits $H$, $k(\beta)$ can be embedded into $H$ \cite[Corollary~3.7,
p.~103]{Mi97}, and we can apply Proposition~10 in \cite[p.~144]{Bou}
to conclude that the image of $\beta$ in $D$ is $c=c_1 + c_2 i + c_3 ij +
c_4 ij$ with reduced trace $\tr(c)=2x_1$ and reduced norm $\nr(c)=pq$.
Hence $2c_1=2x_1$, so $c_1=x_1$ and $x_1 \in T$. Then $(pq)^rx_1 \in S.$
Thus $S \subseteq R_{\pp} \cap R_{\qq}$ and $S$ contains the subgroup $G:=
p^{r+1}R_{\pp} \cap q^{r+1}R_{\qq}$ which has finite index in
$R_{\pp} \cap R_{\qq}$.  This proves the claim.

Let $s_1,\dots,s_{l}$ be coset representatives for $G$ in $R_{\pp}
\cap R_{\qq}$. Then for $x \in k,$
\[
x \in R_{\pp} \cap R_{\qq} \Leftrightarrow (\exists\, s \in
S)(x=s+s_1) \lor \dots \lor (x=s+s_{l}).
\]
This proves that $R_{\pp} \cap R_{\qq}$ is diophantine over $k$.

We can repeat the same argument with $\pp$ and some other finite prime $\ell
\neq \qq$ and conclude that $R_{\pp} \cap R_{\ell}$  is diophantine
over $k$. By weak approximation we have
\[
R_{\pp} = (R_{\pp} \cap R_{\qq})+(R_{\pp} \cap R_{\ell}).
\]
This proves the theorem when the characteristic of $k$ is not $2$.

{\bf Characteristic 2 Case:}
When $k$ has characteristic $2$, we can still find a $4$-dimensional
central division algebra ramified exactly at $\pp$ and $\qq$. We only
have to change the definition of $T$ to 
\[
T:= \{x_3 \in k: (\exists \, x_1,x_2,x_4 \in
k):\nr(x_1+x_2u+x_3v+x_4uv) = pq)\}.
\]
Then we can still show $T \subseteq A_{\pp}$.
For the other direction, given $x_3 \in k$ with $x_3 \in pR_{\pp} \cap q
R_{\qq}$, we look at the equation
\[
X^2 - x_3 X + pq = 0.
\]
Then the proof proceeds exactly as before.
\end{proof}
\section{Integrality at a prime for the perfect closure of global
  fields of characteristic $p>2$}
\label{sectionperf}

{\bf Notation.} In the following $\F_q$ will be the finite field with
$q=p^m$ elements of characteristic $p >2$, $\F_q(t)$ will denote the
field of rational functions over $\F_q$ and $K$ will denote the
perfect closure of $\F_q(t)$, {\it i.e.}\ $K = \F_q(t, t^{1/p},
t^{1/{p^2}}, t^{1/{p^3}}, \cdots)$.
For simplicity of notation we will first prove Theorem~\ref{perfect}
for the rational function field $\F_q(t)$, and then say how the proof
has to be modified for finite extensions $k$ of $\F_q(t)$.
\begin{thm}\label{specialcase}
  Let $K$ be as above.  Let $\pp$ be a prime of $K$. The set $\{ x \in
  K : \ord_{\pp} x \geq 0 \}$ is diophantine over $K$.
\end{thm}
\begin{proof}
  Let $\pp_1$ and $\pp_2$ be two primes of $K$ and let $\ord_{\pp_1}$ and
  $\ord_{\pp_2}$ be the associated additive valuations.
  
  We will show that the set $\{ x \in K : \ord_{\pp_1} x \geq 0 \}$ is
  diophantine over $K$.
  
  The restrictions of $\pp_1$ and $\pp_2$ to $\F_q(t)$ are primes of
  $\F_q(t)$.  For simplicity of notation we will denote these
  restrictions again by $\pp_1$ and $\pp_2$.  From Theorem~\ref{exact}
  and Theorem~\ref{localfield} it follows that we can find a central
  division algebra $D/\F_q(t)$ with $[D:\F_q(t)] =4$ which is ramified
  exactly at the primes $\pp_1$ and $\pp_2$.
\[\mbox{Let }\OO_{D} := \{ z \in D : \ord_{\pp_1}(z)\geq 0
\mbox{ and }\ord_{\pp_2}(z) \geq 0 \},\]
\[\mbox{and }\OO := \{z \in\F_q(t) : \ord_{\pp_1}(z) \geq 0 {\mbox{ and }}
\ord_{\pp_2}(z) \geq 0\}.\]

The ring $\OO$ is an intersection of discrete valuation rings, so
$\OO$ is a Dedekind domain with finitely many primes. By
\cite[Exercise 15, p.~625]{Jac89} $\OO$ is a PID. The ring $\OO_D$ is
a finitely generated torsion-free $\OO$-module. Since $\OO$ is a PID,
it follows that $\OO_{D}$ is a free $\OO$-module of rank 4.

Let $\tr: \OO_{D} \to \OO$ be the reduced trace. Then $\tr(1) =2$,
because $[D : \F_q(t)]=4$. Since $2$ is a unit in $\OO$, the reduced
trace is surjective. Since $\OO_{D}/\OO$ is free, the kernel of the
reduced trace is free of rank $3$, so let $a_2,a_3,a_4$ be a basis for
the kernel. The image of the trace is generated by $\tr(1)$, so
$a_1=1,a_2,a_3,a_4$ are a basis of $\OO_{D}/\OO$. Then
$a_1,\cdots,a_4$ are also a basis for $\OO_D \tensor_\OO \F_q(t)= D$
over $\F_q(t)$.  Let \[S := \{x_1 \in \F_q(t) : (\exists \,
x_2,x_3,x_4 \in \F_q(t)): (\nr(x_1a_1+x_2a_2+x_3a_3+x_4a_4)=1) \}.\]
Then $S \subseteq \OO.$ Let $K= \F_q(t,t^{1/p},t^{1/p^2}, t^{1/p^3},
\cdots)$.

Let $D^{\perf} := D \tensor_{\F_q(t)} K$. Then $D^{\perf}$ is still
ramified at $\pp_1$ and $\pp_2$, because only elements of order
$p^{\ell}$ in $\Br(\F_q(t))$ get killed in the perfection, $D$ has
order 2 in $\Br(\F_q(t))$, and $p \geq 3$.
\[\mbox{Let }\OO^{\perf} := \{z
\in K : \ord_{\pp_1}(z) \geq 0 \mbox{ and } \ord_{\pp_2}(z) \geq 0 \},\]
 \[\mbox{and }\OO_{D^{\perf}}:= \{z \in D^{\perf}: \ord_{\pp_1}(z) 
 \geq 0 \mbox{ and } \ord_{\pp_2}(z) \geq 0\} .\] We will prove that
 $\OO^{\perf}$ is diophantine over $K$. To do this let
\[T:= \{ x_1 \in K: (\exists \, x_2,x_3,x_4 \in K):
(\nr(x_1a_1+x_2a_2+x_3a_3+x_4a_4)=1) \}.
\]
We will prove that $\OO^{\perf}$ is diophantine by showing that there
exist finitely many elements $\alpha_1, \dots,\alpha_r \in K$ such
that
\[
\OO^{\perf} = (T+\alpha_{1}) \cup (T+\alpha_{2}) \cup \cdots \cup
(T+\alpha_{r}).
\]
First we need the following claim:
\\
{\bf Claim:} $\OO_{D^{\perf}}$ is a free $\OO^{\perf}$-module of rank
4 with basis $a_1 \tensor 1, \cdots , a_4 \tensor 1$. Also
$a_1\tensor~1, \cdots, a_4 \tensor 1$ are a basis for $D^{\perf}$
over $K$.
\\
{\bf Proof of Claim:} For each $i \in \N$ let
\begin{align*}
  &D_i := D \tensor_{\F_q(t)} \F_q(t^{1/p^i}), \\&\OO_i:= \{z \in
  \F_q(t^{1/p^i}): \ord_{\pp_1}(z) \geq 0\mbox{ and }\ord_{\pp_2}(z) \geq 0
  \},\mbox{and }\\&\OO_{D_i} := \{z \in \F_q(t^{1/p^i}): \ord_{\pp_1}(z) \geq
  0 \mbox{ and }\ord_{\pp_2}(z) \geq 0\}= \OO_D \tensor_{\OO}
  \OO_i.\end{align*}

Then $\OO_{D_i}$ is a free $\OO_i$-module of rank 4 with basis $a_1
\tensor 1, \cdots , a_4 \tensor 1$ by \cite[Proposition~4.1,
p.~623]{La}.

We have that $\OO_{D^{\perf}} = \OO_D \tensor_{\OO} \OO^{\perf}$, and
hence the same Proposition implies that $\OO_{D^{\perf}}$ is free over
$\OO^{\perf}$ with basis $a_1 \tensor 1, \cdots , a_4 \tensor 1$.
These elements are still linearly independent over the quotient field
of $\OO^{\perf}$, $K$, so they also form a basis for $D^{\perf}$ over
$K$. This proves the claim.

By definition of $T$, we have that $T \subseteq \OO^{\perf}$.  Let
$k_{1}$ and $k_2$ be the residue fields of $\pp_1$ and $\pp_2$,
respectively. The fields $k_1$ and $k_2$ are finite extensions of
$\F_q$.  For $x_1 \in \OO^{\perf}$ we have:
\begin{eqnarray}
\nonumber & &{x_1^2-1} \mod \pp_i \notin (k_i)^2 \mbox{ for }
i=1,2\\
\label{Eq:3}& \Rightarrow & x_1^2 -1 \notin {(K_v^*)}^2 \mbox{ locally at }
 v=\pp_1,\pp_2\\
\label{Eq:4}&\Leftrightarrow & \begin{cases} X^2 - 2x_1X +1 
\mbox{ is irreducible over } K_v
\mbox{ for }v=\pp_1,\pp_2 \\\mbox{ or }x_1 = \pm 1 \end{cases}\\
\label{Eq:1}& \Leftrightarrow & x_1=\pm 1 \mbox{ or }
(\exists \, \alpha \in D^{\perf} \mbox{ s.t. }K(\alpha) \mbox{ splits }
 D^{\perf}, \\\nonumber&&\mbox{ and } \alpha^2-2\alpha x_1+
1=0)\\
\label{Eq:2}& \Leftrightarrow &  x_1 = \pm 1 \mbox{ or }(\exists\, \alpha \in
D^{\perf}\mbox{ s.t. }\tr(\alpha) = 2 x_1,
   \nr (\alpha) =1,\\\nonumber&&\mbox{ and }[K(\alpha) :K] = 2)\\
\nonumber & \Leftrightarrow & \exists \,\alpha \in D^{\perf} \mbox{
  s.t. }\tr(\alpha) = 2 x_1, \mbox{ and } \nr (\alpha) =1 \\ 
\nonumber & \Leftrightarrow & x_1 \in T.
\end{eqnarray}
The equivalence of (\ref{Eq:3}) and (\ref{Eq:4}) comes from solving
the equation $X^2 - 2x_1X +1$ using the quadratic formula.  The
equivalence of (\ref{Eq:1}) and (\ref{Eq:2}) follows from the fact
that every degree $2$ field extension $K(\alpha) \subseteq D^{\perf}$
splits the 4-dimensional division algebra $D^{\perf}$.

There exists an $a_1 \in k_1$ such that $(a_1^2-1) \notin (k_1)^2$: If
$a_1^2-1$ were a square for every $a_1 \in k_1$, then we would have
$a_1^2 -1 = b^2 $, so $a_1^2 -2 = b^2 -1 = c^2$ is a square, so
repeating this $p$ times for every square we could show that the
number of squares in $k_1$ is divisible by $p$.  But $k_1=\F_{p^n}$
for some $n>0$ and the number of squares in $\F_{p^n}$ is $(p^n+1)/2$
which is not divisible
by $p$. \\
The same argument shows that there exists an element $a_2 \in k_2$
such that $(a_2^2-1) \notin (k_2)^2$.

Let $a_1 \in k_1$ and $a_2 \in k_2$ be such elements.  By the
approximation theorem there exists an element $a \in \OO^{\perf}$ such
that $a \equiv a_1 \mod \pp_1$ and $a \equiv a_2
\mod \pp_2$.  From the above equivalences it follows
that $a \in T$.  The approximation theorem implies that for each $i
\in k_1$, $j \in k_2$ we can find an element $\alpha_{i,j} \in
\OO^{\perf}$ with the property that $\alpha_{i,j}\equiv i \mod \pp_1$
and $\alpha_{i,j} \equiv j \mod \pp_2$.
\\
{\bf Claim:} \[\OO^{\perf} = \bigcup_{i \in k_1,j \in
  k_2}(T+\alpha_{i,j}).\] {\bf Proof of Claim:} The set $T$ contains
all elements \[\{ x \in K: x \equiv a_1 \mod \pp_1 \mbox{ and } x
\equiv a_2 \mod \pp_2 \}.\] If $y \in \OO^{\perf}$, then for some $i
\in k_1,j \in k_2$, $y \equiv i \mod \pp_1$ and $y \equiv j \mod
\pp_2$, so then $y - \alpha_{(i-a_1),(j-a_2)} \in T$. This proves the
claim.

The claim implies that $\OO^{\perf}$ is diophantine over $K$. The same
argument with $\pp_2$ replaced by some other prime $\pp_3$ shows that
the set $\tilde {\OO}^{\perf }=\{ z \in K: \ord_{\pp_1}(z) \geq 0$ and $
\ord_{\pp_3} \geq 0\}$ is diophantine over $K$. Then by weak
approximation $\{x \in K: \ord_{\pp_1}(x) \geq 0 \} = \OO^{\perf } +
\tilde{\OO}^{\perf}$.
\end{proof}
\begin{lem}
  Let $k$ be any global field of characteristic $p>0$ such that $k$ is
  a finite extension of $\F_q(t)$ for some $q=p^n$. The perfect
  closure of $k$ is $k^{\perf}:= k (t^{1/p},
  t^{1/{p^2}},t^{1/{p^3}},\dots)$.
\end{lem}
\begin{proof}
  Clearly $k^{\perf}$ is contained in the perfect closure of $k$. The
  field $k^{\perf}$ is a finite extension of $K = \F_q(t, t^{1/p},
  t^{1/{p^2}}, t^{1/{p^3}}, \cdots)$. Since $K$ is perfect, and finite
  extensions of perfect fields are perfect, $k^{\perf}$ is perfect as
  well, so it must be equal to the perfect closure of $k$.
\end{proof}
Now we can state the general theorem:
\begin{thm}
  Let $k$ be a global field of characteristic $p>2$, and $k^{\perf}$
  its perfect closure. Let $\pp$ be a prime of $k^{\perf}$.  The set
  $\{ x \in k^{\perf}: \ord_{\pp} x \geq 0 \}$ is diophantine over
  $k^{\perf}$.
\end{thm}
\begin{proof}
  We can repeat the proof of Theorem~\ref{specialcase} with $\F_q(t)$
  replaced by $k$. Everything works exactly as before, because the
  exact sequence of Theorem~\ref{exact} works for all global fields
  $k$.
\end{proof}
\noindent{\bf Acknowledgments.} I thank Bjorn Poonen
for suggesting the method of proof of Theorem~\ref{orderproof} and for
his comments regarding Section~\ref{sectionperf}.

\end{document}